\documentclass[10pt]{elsart}
\usepackage{amssymb,amsmath,amstext,graphicx,subfigure,warmread,times}
\usepackage[all,import]{xy}

\newcommand{\norm}[1]{\ensuremath{\left\| #1 \right\|}}
\newcommand{\bracket}[1]{\ensuremath{\left[ #1 \right]}}
\newcommand{\braces}[1]{\ensuremath{\left\{ #1 \right\}}}
\newcommand{\parenth}[1]{\ensuremath{\left( #1 \right)}}
\newcommand{\refeqn}[1]{(\ref{eqn:#1})}
\newcommand{\reffig}[1]{Fig. \ref{fig:#1}}
\newcommand{\tr}[1]{\mbox{tr}\ensuremath{\bracket{#1}}}
\newcommand{\deriv}[2]{\ensuremath{\frac{\partial #1}{\partial #2}}}
\newcommand{\SO}{\ensuremath{\mathrm{SO(3)}}}
\newcommand{\T}{\ensuremath{\mathrm{T}}}
\newcommand{\so}{\ensuremath{\mathfrak{so}(3)}}

\renewcommand{\Re}{\ensuremath{\mathbb{R}}}
\renewcommand{\S}{\ensuremath{\mathbb{S}}}

\newcommand{\wh}{\widehat}
\newcommand{\wt}{\widetilde}
\newcommand{\la}{\label}
\newcommand{\be}{\begin{equation}}
\newcommand{\ee}{\end{equation}}
\newcommand{\bea}{\begin{eqnarray}}
\newcommand{\eea}{\end{eqnarray}}
\newcommand{\beas}{\begin{eqnarray*}}
\newcommand{\eeas}{\end{eqnarray*}}
\newcommand{\nn}{\nonumber}

\newcommand{\Tp}{^{\mbox{\rm\small T}}}
\newtheorem{theorem}{Theorem}[section]
\newtheorem{proposition}{Proposition}[section]
\newtheorem{lemma}{Lemma}[section]

\begin{document}

\begin{frontmatter}

\title{Global Optimal Attitude Estimation using Uncertainty Ellipsoids}
\author{Amit K. Sanyal\corauthref{cor1}\thanksref{la1}, Taeyoung Lee\thanksref{la2},
Melvin Leok\thanksref{la3}, N. Harris McClamroch\thanksref{la2}}
\ead{sanyal@asu.edu, \{tylee, mleok, nhm\}@umich.edu}
\corauth[cor1]{Corresponding author. Tel: 1-480-727-8783.}
\address[la1]{Department of Mechanical and Aerospace Engineering,
Arizona State University, Tempe, AZ 85287-6106}
\address[la2]{Department of Aerospace Engineering, University of Michigan,
Ann Arbor, MI 48109-2140}
\address[la3]{Department of Mathematics, University of Michigan,
Ann Arbor, MI 48109-1043}

\begin{abstract}
A deterministic attitude estimation problem for a rigid
body in a potential field, with bounded attitude and angular velocity measurement
errors is considered. An attitude estimation algorithm that globally minimizes
the attitude estimation error is obtained.
Assuming that the initial attitude, the initial angular velocity and measurement
noise lie within given ellipsoidal bounds, an uncertainty ellipsoid that bounds
the attitude and the angular velocity of the rigid body is obtained. The center
of the uncertainty ellipsoid provides point estimates, and the size of the
uncertainty ellipsoid measures the accuracy of the estimates. The point
estimates and the uncertainty ellipsoids are propagated using a Lie group
variational integrator and its linearization, respectively. The attitude
and angular velocity estimates are optimal in the sense that the sizes of the
uncertainty ellipsoids are minimized.
\end{abstract}

\end{frontmatter}

\section{Introduction}
\vspace*{-3mm}

Attitude estimation is often a prerequisite for controlling aerospace and
underwater vehicles, mobile robots, and other mechanical systems moving in
space. In this paper, we study the attitude estimation problem for the
uncontrolled dynamics of a rigid body in an attitude-dependent force potential
(like uniform gravity). The estimation scheme we present has the following
important features: (1) the attitude is globally represented without using
any coordinate system, (2) the filter obtained is not a Kalman or extended
Kalman filter, and (3) the attitude and angular velocity measurement errors
are assumed to be bounded, with known ellipsoidal uncertainty bounds. The
static attitude estimation scheme presented here is based on~\cite{pro:san2006}.
The attitude is represented globally using proper orthogonal matrices and
the exponential map on the set of $3\times 3$ skew-symmetric matrices. Such
a global representation has been recently used for partial attitude estimation
with a linear dynamics model in~\cite{jo:rehu}. The estimation scheme
presented here is deterministic, based on known measurement uncertainty bounds
propagated by the dynamic flow.

The attitude determination problem for a rigid body from vector measurements was
first posed in \cite{jo:wahba}. A sample of the literature in spacecraft attitude
estimation can be found in \cite{jo:bosh,jo:cram,jo:mark,jo:shus1,jo:shus2}.
Applications of attitude estimation to unmanned vehicles and robots can be found
in \cite{jo:rehu,jo:badu,pro:rosb,pro:vafo}. Most existing attitude estimation
schemes use generalized coordinate representations of the attitude. As is well
known, minimal coordinate representations of the rotation group, like Euler angles,
Rodrigues parameters, and modified Rodrigues parameters (see \cite{jo:cram2}),
usually lead to geometric or kinematic singularities. Non-minimal coordinate
representations, like the unit quaternions used in the quaternion estimation
(QUEST) algorithm \cite{jo:shus1} and its several variants
\cite{jo:bosh,jo:shus2,jo:psia}, have their own associated problems. Besides
the extra constraint of unit norm that one needs to impose, the quaternion
vector itself can be defined in one of two ways, depending on the sense of
rotation used to define the principal angle.

A brief outline of this paper is given here. In Section 2, the attitude
determination problem for vector measurements with measurement noise is
introduced, and a global attitude determination algorithm which minimizes
the attitude estimation error is presented. The attitude dynamics and dynamic
estimation problem is formulated for a rigid body in an attitude-dependent
potential. Section 3 presents the attitude estimation scheme assuming that
both attitude and angular velocity measurements are available simultaneously.
Sufficient conditions for convergence of the estimates are given. This
attitude estimation scheme has also been extended and applied to the case
where only attitude but no angular velocity measurements are available,
and recently reported in~\cite{pro:CDC06}. Section 4 presents concluding
remarks and observations.

\section{Attitude estimation from vector observations}
\vspace*{-3mm}

Attitude of a rigid body is defined as the orientation of a body
fixed frame with respect to an inertial reference frame; it is
represented by a rotation matrix that is a $3\times 3$ orthogonal
matrix with determinant $+1$. Rotation matrices have a group structure
denoted by $\SO$, and its action on $\Re^3$ takes a vector represented in body fixed
frame into its representation in the reference frame by matrix
multiplication.
\vspace*{-2mm}
\subsection{Attitude determination procedure}
\vspace*{-2mm}

We assume that there are $m$ fixed points in the spatial reference
frame, no two of which are co-linear, that are measured in the body
frame. We denote the known direction of the $i$th point in the spatial
reference frame as $e^i\in\S^2$, and the corresponding vector
represented in the body fixed frame as $b^i\in\S^2$. Since we only
measure directions, we normalize $e^i$ and $b^i$ so that they have
unit lengths. The $e^i$ and $b^i$ are related by a rotation matrix
$C\in\SO$ that defines the attitude of the rigid body;
\begin{align*}
e^i = C b^i,
\end{align*}
for all $i\in\braces{1,2,\cdots,m}$. We assume that $e^i$ is known
accurately and we assume that $b^i$ is measured in the body fixed
frame. Let the measured direction vector be $\tilde{b}^i\in\S^2$,
which contains measurement errors, and let an estimate of the
rotation matrix be $\wh{C}\in\SO$. The vector estimation errors are
given by the
\begin{align*}
e^i-\wh{C}\tilde{b}^i,\ i=1,\ldots,m.
\end{align*}

The attitude determination problem consists of finding
$\wh{C}\in\SO$ such that the weighted $2$ norm of these errors is
minimized:
\begin{align}
\min_{\wh{C}} \mathcal J & =\frac{1}{2}\sum_{i=1}^m w_i
(e^i-\wh{C}\tilde b^i)\Tp (e^i-\wh{C}\tilde b^i),\nonumber\\
& = \frac{1}{2} \tr {(E-\wh{C}\tilde B)\Tp W(E-\wh{C}\tilde B)},\label{eqn:wahba}\\
&\text{subject to } \wh{C}\in\SO,\nonumber
\end{align}
where $E=\bracket{e^1,e^2,\cdots,e^m}\in\Re^{3\times m}$, $\tilde
B=\bracket{\tilde b^1,\tilde b^2,\cdots,\tilde b^m}\in\Re^{3\times
m}$, and $W=\mathrm{diag}\bracket{w^1,w^2,\cdots,w^m}\in\Re^{m\times
m}$ has a weighting factor for each measurement. We assume that $m\ge 3$
in this paper. If $m=2$, we can take the cross-product of the two
measured unit vectors $\tilde b^1\times\tilde b^2=\tilde b^3$ and treat
that as a third measured direction; the corresponding unit vector in
the inertial frame is $e^3= e^1\times e^2$.

The problem \refeqn{wahba} is known as Wahba's problem~\cite{jo:wahba}.
The original solution of Wahba's problem is given in~\cite{jo:solwahba},
and a solution expressed in terms of quaternions, known as the QUEST
algorithm, is presented in~\cite{jo:shus1}. A solution without using
generalized attitude coordinates is given in~\cite{pro:san2006}. A
necessary condition for optimality of \refeqn{wahba} is given
in~\cite{pro:san2006} as
\begin{align}
L\Tp \wh{C}=\wh{C}\Tp L,\label{eqn:san2006}
\end{align}
where $L=EW\wt B\Tp \in\Re^{3\times 3}$.

The following result, proved in~\cite{pro:san2006}, gives an unique
estimate $\wh{C}\in\SO$ of the attitude matrix that satisfies
\refeqn{san2006} and solves the attitude determination problem
\refeqn{wahba}.
\begin{theorem}
The unique minimizing solution to the attitude determination problem
\refeqn{wahba} is given by
\be \wh{C}= SL,\;\ S=Q\sqrt{(RR\Tp)^{-1}}Q\Tp, \la{soln} \ee
where
\be L=QR,\;\ Q\in\SO, \la{LQR} \ee
and $R$ is upper triangular and invertible; this is the QR decomposition
of $L$. The symmetric positive definite (principal) square root is used
in (\ref{soln}).
\la{atde}
\end{theorem}

\vspace*{-2mm}
\subsection{State bounding estimation}
\vspace*{-2mm}

Here we discuss the general idea of deterministic state bounding
estimation, using ellipsoidal sets to describe state uncertainty
and measurement noise. A stochastic state estimator requires
probabilistic models for the state uncertainty and the noise.
However, statistical properties of the uncertainty and the noise are
often not available. We usually make statistical assumptions on
disturbance and noise in order to make the estimation problem
mathematically tractable. In many practical situations such idealized
assumptions are not appropriate, and this may cause poor estimation
performance~\cite{jo:TheSkaSou.IEEETSP94}.

An alternative deterministic approach is to specify bounds on the
uncertainty and the measurement noise without any assumption on its
distribution. Noise bounds are available in many cases, and
deterministic estimation is robust to the noise distribution. An
efficient but flexible way to describe the bounds is using
ellipsoidal sets, referred to as uncertainty ellipsoids.

The deterministic estimation process using uncertainty ellipsoids
has the same structure as the Kalman filter. We assume that the
initial state lies in a prescribed uncertainty ellipsoid. We
propagate the initial uncertainty ellipsoid through time using the
equations of motion. This defines a prediction step. After measuring
the state, we find a bound on the state assuming that measurement
error is bounded. The bound is described by a measurement uncertainty
ellipsoid. We then obtain a minimal ellipsoid that contains the intersection of the
predicted uncertainty ellipsoid and the measured uncertainty ellipsoid.
This procedure is repeated whenever new measurements are available.

\renewcommand{\xyWARMinclude}[1]{\includegraphics[width=0.45\textwidth]{#1}}
\begin{figure*}[t]
    \centerline{\subfigure[Propagation of uncertainty ellipsoid]{
    $$\begin{xy}
    \xyWARMprocessEPS{tube}{eps}%%
    \xyMarkedImport{}%%
    \xyMarkedMathPoints{1-5}
    \end{xy}$$}
    \hspace*{2cm}
    \renewcommand{\xyWARMinclude}[1]{\includegraphics[width=0.3\textwidth]{#1}}
    \subfigure[Filtering procedure]{
        $$\begin{xy}
    \xyWARMprocessEPS{tubesection}{eps}%%
    \xyMarkedImport{}%%
    \xyMarkedMathPoints{1-3}
    \end{xy}$$}}
%    }%%
    \caption{Uncertainty ellipsoids}\label{fig:ue}
\end{figure*}
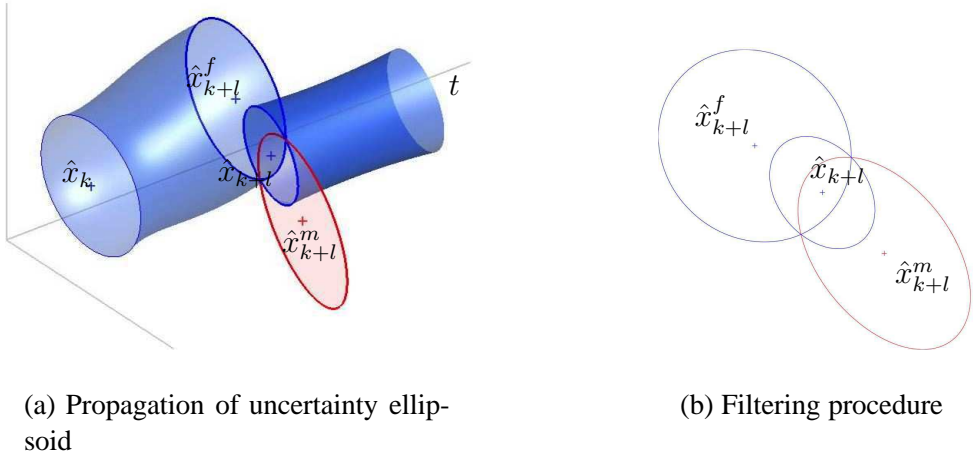

This deterministic estimation procedure is illustrated in \reffig{ue},
where the left figure shows evolution of an uncertainty ellipsoid in
time, and the right figure shows a cross section at a fixed time when
the state is measured. Suppose that the time interval between two sets
of measurements is divided into $l$ equal time steps for discrete
integration, and the subscript $k$ denotes the $k$th discrete integration
time step. At the previous measurement instant, corresponding to the
$k$th time step, the state is bounded by an uncertainty ellipsoid
centered at the estimated state $\hat{x}_k$. This initial ellipsoid is
propagated through time. Depending on the dynamics of the system, the
size and the shape of the tube are changed. The new set of measurements
are taken at the $(k+l)$th time step. At the $(k+l)$th time step, the
predicted uncertainty ellipsoid is centered at $\hat{x}_{k+l}^f$. Another
ellipsoidal bound on the state is obtained from the measurements. The
measured uncertainty
ellipsoid is centered at $\hat{x}_{k+l}^m$. The state lies in the
intersection of the two ellipsoids. In the estimation procedure, we
find a new ellipsoid that contains the intersection, which is shown in
the right figure. The center of the new ellipsoid, $\hat{x}_{k+l}$ is
considered as a point estimate at time step $k+l$, and the magnitude
of the new uncertainty ellipsoid measures the accuracy of the estimation.
If the size of the uncertainty ellipsoid is small, then we can conclude
that the estimated state is accurate. The deterministic estimation is
optimal in the sense that the size of the new ellipsoid is minimized.

The deterministic estimation process is based on the state estimation techniques
developed in~\cite{jo:Sc1968}. Optimal deterministic state or parameter
estimation is considered in~\cite{jo:MaNo1996,jo:DuWaPo2001,jo:PoNaDuWa2004},
where an analytic solution for the minimum ellipsoid that contains a union
or an intersection of ellipsoids is given. Parameter estimation in the
presence of bounded noise is dealt with in~\cite{jo:BaNaTe96,jo:BeGa2004}.

\vspace*{-2mm}
\subsection{Attitude Estimation Problem formulation}
\vspace*{-2mm}

\subsubsection{Equations of motion}
We consider estimation of the attitude dynamics of a rigid body
in the presence of an attitude dependent potential. We assume that
the potential $U(\cdot):\SO\mapsto\Re$ is determined by the attitude
of the rigid body, $C\in\SO$. A spacecraft on a circular orbit
including gravity gradient effects~\cite{pro:acc06}, or a 3D
pendulum~\cite{pro:cca05} can be modeled in this way. The continuous
equations of motion are given by
\begin{gather}
J\dot\omega + \omega\times J\omega = M,\\
\dot{C} = C S(\omega),\label{eqn:Rdot}
\end{gather}
where $J\in\Re^{3\times 3}$ is the moment of inertia matrix of the
rigid body, $\omega\in\Re^3$ is the angular velocity of the body
expressed in the body fixed frame, and $S(\cdot):\Re^3\mapsto \so$
is a skew mapping defined such that $S(x)y=x\times y$ for all
$x,y\in\Re^3$. $M\in\Re^3$ is the moment due to the potential. The
moment is determined by the relationship, $S(M)=\deriv{U}{C}\Tp
C-C\Tp\deriv{U}{C}$, or more explicitly,
\begin{gather}
M=r_1\times v_{r_1} + r_2\times v_{r_2} +r_3\times v_{r_3},
\end{gather}
where $r_i,v_{r_i}\in\Re^{1\times 3}$ are the $i$th row vectors of
$C$ and $\deriv{U}{C}$, respectively. The detailed description of
this rigid body model and the derivation of the above equations can
be found in~\cite{pro:cca05}.

General numerical integration methods, including the popular
Runge-Kutta schemes, typically preserve neither first integrals nor
the characteristics of the configuration space, $\SO$. In particular,
the orthogonal structure of the rotation matrices is not preserved
numerically. It is often proposed to parameterize \refeqn{Rdot} by
Euler angles or quaternions instead of integrating \refeqn{Rdot} directly.
However, Euler angles have singularities. The numerical simulation
process has to be monitored and switching between Euler angle charts
is necessary in order to avoid singularities. Quaternions are free of
singularities, but the quaternion representing the attitude is required
to have unit length. The matrix corresponding to a quaternion which is
not of unit length is not orthogonal, and hence does not represent a rotation.

To resolve these problems, a Lie group variational integrator for
the attitude dynamics of a rigid body is proposed in~\cite{pro:cca05}.
This Lie group variational integrator is described by the discrete
time equations.
\begin{gather}
h S(J\omega_k+\frac{h}{2} M_k) = F_k J_d - J_dF_k\Tp,\label{eqn:findf0}\\
C_{k+1} = C_k F_k,\label{eqn:updateR0}\\
J\omega_{k+1} = F_k\Tp J\omega_k +\frac{h}{2} F_k\Tp M_k
+\frac{h}{2}M_{k+1},\label{eqn:updatew0}
\end{gather}
where $J_d\in\Re^3$ is a nonstandard moment of inertia matrix
defined by $J_d=\frac{1}{2}\tr{J}I_{3\times 3}-J$, and $F_k\in\SO$
is the relative attitude over an integration step. The constant
$h\in\Re^+$ is the integration step size. This integrator yields
a map $(C_k,\omega_k)\mapsto(C_{k+1},\omega_{k+1})$ by solving
\refeqn{findf0} to obtain $F_k\in\SO$ and substituting it into
\refeqn{updateR0} and \refeqn{updatew0} to obtain $C_{k+1}$ and
$\omega_{k+1}$. Numerically, we ensure that $F_k$ remains on $\SO$ by requiring that $F_k=\exp(S(f_k))$, where $f_k\in\Re^3$. This allows us to express the discrete equations in terms of $f_k\in\Re^3$ as opposed to $F_k\in\SO$.

Since this integrator does not use a local parameterization, the
attitude is defined globally without singularities. It preserves the
orthogonal structure of $\SO$ because the rotation matrix is updated
by a multiplication of two rotation matrices in \refeqn{updateR0},
which is a group operation of \SO. This integrator is obtained from
a discrete variational principle, and it exhibits the characteristic
symplectic and momentum preservation properties, and good energy
behavior characteristic of variational integrators. We use
\refeqn{findf0}, \refeqn{updateR0}, and \refeqn{updatew0} in the
following development of the attitude estimator.

\vspace*{-2mm}
\subsubsection{Uncertainty Ellipsoid}
\vspace*{-2mm}

An uncertainty ellipsoid in $\Re^n$ is defined as
\begin{align}
    \mathcal{E}_{\Re^n}(\hat x,P)=\braces{x\in\Re^n \,\Big|\,
    (x-\hat{x})\Tp P^{-1}(x-\hat{x})\leq
    1},
\end{align}
where $\hat{x}\in\Re^n$, and $P\in\Re^{n\times n}$ is a symmetric
positive definite matrix. We call $\hat x$ the center of the
uncertainty ellipsoid, and we call $P$ the uncertainty matrix that
determines the size and the shape of the uncertainty ellipsoid. The
size of an uncertainty ellipsoid is measured by $\tr{P}$. It equals
the sum of the squares of the semi principal axes of the ellipsoid.

The configuration space of the attitude dynamics is $\SO$, so the
state evolves in the 6 dimensional tangent bundle, $\T\SO$. Thus the
corresponding uncertainty ellipsoid is a submanifold of $\T\SO$. An
uncertainty ellipsoid centered at $(\hat{C},\hat\omega)$ is induced
from an uncertainty ellipsoid in $\Re^6$, using the Lie algebra $\so$;
\begin{align}
    \mathcal{E}(\hat{C},\hat{\omega},P) & = \braces{C\in\SO,\,\omega\in\Re^3 \,\Big|\,
    \begin{bmatrix}\zeta\\\delta\omega\end{bmatrix}\in\mathcal{E}_{\Re^6}(0_{6},P)},\label{eqn:ueso}
\end{align}
where $S(\zeta)=\mathrm{logm} \parenth{\hat{C}\Tp C}\in\so$,
$\delta\omega=\omega-\hat{\omega}\in\Re^3$, and $P\in\Re^{6\times
6}$ is a symmetric positive definite matrix. Equivalently, an
element $(C,\omega)\in\mathcal{E}(\hat{C},\hat{\omega},P)$ can be
written as
\begin{align*}
    C & = \hat{C} \exp{\big(S(\zeta)\big)},\\
    \omega & = \hat{\omega} + \delta \omega,
\end{align*}
for some $x=[\zeta\Tp,\,\delta\omega\Tp]\Tp\in\Re^6$ satisfying
$x\Tp P^{-1}x\leq 1$.

\vspace*{-2mm}
\subsubsection{Uncertainty model}
\vspace*{-2mm}
We describe the measurement error models for the measured direction
vectors and the angular velocity. The direction vector $b^i\in\S^2$
is measured in the body fixed frame, and let $\tilde b^i\in\S^2$
denote the measured direction. Since we only measure directions, we
normalize $b^i$ and $\tilde b^i$ so that they have unit lengths.
Therefore it is inappropriate to express the measurement error by a
vector difference. The measurement error is modeled by rotation of
the measured direction;
\begin{align}
{b}^i& = \exp{\big(S(\nu^i)\big)} \tilde b^i\nonumber\\
& \simeq \tilde b^i + S(\nu^i)\tilde b^i,\label{eqn:bi}
\end{align}
where $\nu^i\in\Re^3$ is the measurement error, which represents the
Euler axis of rotation vector from $\tilde b^i$ to $b^i$, and
$\norm{\nu^i}$ is the corresponding rotation angle in radians. We
assume that the measurement error is small to obtain the second
equality.

The angular velocity measurement errors are modeled as
\begin{align}
\omega=\tilde\omega + \upsilon,\label{eqn:Omega}
\end{align}
where $\tilde\omega \in\Re^3$ is the measured angular velocity,
and $\upsilon\in\Re^3$ is an additive measurement error.

We assume that the initial conditions and the measurement error are
bounded by prescribed uncertainty ellipsoids.
\begin{gather}
    (C_0,\omega_0)\in\mathcal{E}(\hat{C}_0,\hat{\omega}_0,P_0),\label{eqn:P0}\\
    \nu^i\in\mathcal{E}_{\Re^3}(0,S^i),\label{eqn:Sk}\\
    \upsilon\in\mathcal{E}_{\Re^3}(0,T)\label{eqn:Tk},
\end{gather}
where $P_0\in\Re^{6\times 6}$, and $S^i, T \in\Re^{3\times 3}$ are
symmetric positive definite matrices that define the shape and the
size of the uncertainty ellipsoids.

\section{Attitude Estimation with Angular Velocity Measurements}
\la{atvfilt}\vspace*{-3mm}

In this section, we develop a deterministic estimator for the
attitude and the angular velocity of a rigid body assuming that both
the attitude measurement and the angular velocity measurements are
available. The estimation process consists of three stages; flow
update, measurement update, and filtering. The flow update predicts
the uncertainty ellipsoid in the future. The measurement update
finds an uncertainty ellipsoid in the state space using the
measurements and the measurement error model. The filtering stage obtains
a new uncertainty ellipsoid compatible with the predicted uncertainty
ellipsoid and the measured uncertainty ellipsoid.

The superscript $i$ denotes the $i$th directional measurement. The
superscript $f$ denotes variables related to the flow update, while
the superscript $m$ denotes variables related to the measurement
update. The notation $\tilde\cdot$ denotes a measured variable, while
$\hat\cdot$ denotes an estimated variable.

\vspace*{-2mm}
\subsection{Flow update}
\vspace*{-2mm}
Suppose that the attitude and the angular momentum at the $k$th step,
which corresponds to the previous measurement instant, lie in a given
uncertainty ellipsoid:
\begin{align*}
    (C_k,\omega_k)\in\mathcal{E}(\hat{C}_k,\hat{\omega}_k,P_k).
\end{align*}

The flow update gives us the center and the uncertainty matrix
that define the uncertainty ellipsoid at the $(k+l)$th step (the
current measurement instant) using the
given uncertainty ellipsoid at the $k$th step. Since the attitude
dynamics of a rigid body is nonlinear, the boundary of the state at
the $(k+l)$th step is not an ellipsoid in general. We assume that the
given uncertainty ellipsoid at the $k$th step is sufficiently small
that the states in the uncertainty ellipsoids can be approximated by
linearized equations of motion. Then we can guarantee that the
boundary of the state at the $(k+l)$th step is an ellipsoid, and we
can compute the center and the uncertainty matrix at the $(k+l)$th
step separately.

\textit{Center:} For the given center, $(\hat{C}_k,\hat{\omega}_k)$,
the center of the uncertainty ellipsoid due to flow propagation is
denoted $(\hat{C}_{k+1}^{f},\hat{\omega}_{k+1}^{f})$. This center is
obtained from the discrete equations of motion, \refeqn{findf0},
\refeqn{updateR0}, and \refeqn{updatew0} applied to $(\hat{C}_k,
\hat{\omega}_k)$:
\begin{gather}
h S(J\hat{\omega}_k+\frac{h}{2} \hat{M}_k) = \hat{F}_k J_d - J_d
\hat{F}_k\Tp,\label{eqn:findf}\\
\hat{C}_{k+1}^{f} = \hat{C}_k \hat{F}_k,\label{eqn:updateR}\\
J\hat{\omega}_{k+1}^{f} = \hat{F}_k\Tp J
\hat{\omega}_k+\frac{h}{2}\hat{F}_k\Tp \hat{M}_k +\frac{h}{2}
\hat{M}_{k+1}.\label{eqn:updatePi}
\end{gather}
This integrator yields a map $(\hat C_k,\hat\omega_k)\mapsto(\hat
C^f_{k+1},\hat\omega_{k+1}^f)$, and this process can be repeated to
find the center at the $(k+l)$th step, $(\hat
C^f_{k+l},\hat\omega_{k+l}^f)$.

\textit{Uncertainty matrix:} We assume that an uncertainty ellipsoid
contains small perturbations from its center. Then the uncertainty
matrix is obtained by linearizing the above discrete equations of
motion. At the $(k+1)$th step, the uncertainty ellipsoid is
represented by perturbations from the center $(\hat C^f_{k+l},
\hat\omega_{k+l}^f)$ as
\begin{align*}
    C_{k+1}&=\hat{C}_{k+1}^{f} \exp{\Big(S(\zeta_{k+1}^{f})\Big)},\\
    \omega_{k+1}&=\hat{\omega}_{k+1}^{f}+\delta\omega_{k+1}^{f},
\end{align*}
for some $\zeta_{k+1}^{f},\delta\omega_{k+1}^{f}\in\Re^3$. The
uncertainty matrix at the $(k+1)$th step is obtained by finding a
bound on $\zeta_{k+1}^{f},\delta\omega_{k+1}^{f}\in\Re^3$. Assume
that the uncertainty ellipsoid at the $k$th step is sufficiently
small. Then, $\zeta_{k+1}^{f},\delta\omega_{k+1}^{f}$ are
represented by the following linear equations (using the results
presented in~\cite{pro:acc06}):
\begin{align*}
x_{k+1}^{f} & = A_k^f x_k,
\end{align*}
where $x_k=[\zeta_k\Tp,\delta\omega_k\Tp]\Tp\in\Re^6$, and
$A_k^f\in\Re^{6\times 6}$ can be suitably defined. Since
$(C_k,\omega_k)\in\mathcal{E}(\hat{C}_k,\hat{\omega}_k,P_k)$,
$x_k\in\mathcal{E}_{\Re^6}(0,P_k)$ by the definition of the
uncertainty ellipsoid given in \refeqn{ueso}. Then we can show that
$A_k^f x_k$ lies in the following uncertainty ellipsoid.
\begin{align*}
A_k^f x_k&\in\mathcal{E}_{\Re^6}\!\parenth{0,A_k^f P_k
\parenth{A_k^f}\Tp}.
\end{align*}
Thus, the uncertainty matrix at the $(k+1)$th step is given by
\begin{align}
P_{k+1}^f & = A_k^f P_k \parenth{A_k^f}\Tp.\label{eqn:Pkpf}
\end{align}
The above equation can be applied repeatedly to find the uncertainty
matrix at the $(k+l)$th step.

We have obtained expressions to predict the center and the
uncertainty matrix in the future from the current uncertainty
ellipsoid using the discrete flow. In summary, the uncertainty
ellipsoid at the ($k+l$)th step is computed using \refeqn{findf},
\refeqn{updateR}, \refeqn{updatePi}, and \refeqn{Pkpf} as:
\begin{align}\label{eqn:flow}
    (C_{k+l},\omega_{k+l})\in\mathcal{E}(\hat{C}_{k+l}^{f},
\hat{\omega}_{k+l}^{f},P_{k+l}^f),\;\ P_{k+l}^f & = A^f
P_k (A_k^f)\Tp,
\end{align}
where $A^f= A^f_{k+l-1}A^f_{k+l-2}\cdots A^f_k\in\Re^{6\times 6}$.

\vspace*{-2mm}
\subsection{Measurement update}
\vspace*{-2mm}

We assume that the attitude and the angular velocity of a rigid body
are measured simultaneously. The measured attitude and the
measured angular velocity have uncertainties since the measurements
contain measurement errors. However, we can find bounds for the actual
state because the measurement errors are bounded by known uncertainty
ellipsoids given by \refeqn{Sk} and \refeqn{Tk}. The measurement
update stage finds an uncertainty ellipsoid in the state space using
the measurements and the measurement error models. The measured attitude
and the measured angular velocity are the center of the measured
uncertainty ellipsoid, and the measurement error models are used to find
the uncertainty matrix.

\textit{Center:} The center of the uncertainty ellipsoid,
$(\hat{C}_{k+l}^{m},\hat\omega_{k+l}^m)$ is obtained by
measurements. The attitude is determined by measuring the directions
to the known points in the inertial frame. Let the measured
directions to the known points be $\tilde B_{k+l}=\bracket{\tilde
b^1,\tilde b^2,\cdots,\tilde b^m}\in\Re^{3\times m}$. Then, the
attitude $\hat{C}_{k+l}^{m}$ satisfies the following necessary
condition given in \refeqn{san2006}:
\begin{gather}
\parenth{\hat{C}_{k+l}^m}\Tp\tilde L_{k+l}-\tilde L_{k+l}\Tp
\hat{C}_{k+l}^m=0,\label{eqn:meaR}
\end{gather}
where $\tilde L_{k+l}=E_{k+l} W_{k+l}
\tilde{B}_{k+l}\Tp\in\Re^{3\times 3}$. The solution of \refeqn{meaR}
is obtained by a QR factorization of $\tilde L_{k+l}$ as given in
Theorem \ref{atde}
\begin{align}
\hat{C}_{k+l}^m=\parenth{Q_q\sqrt{(Q_rQ_r\Tp)^{-1}}Q_q\Tp} \tilde
L_{k+l},\label{eqn:meaRex}
\end{align}
where $Q_q\in\SO$ is an orthogonal matrix and $Q_r\in\Re^{3\times
3}$ is a upper triangular matrix satisfying $\tilde L_{k+l}=Q_qQ_r$.

The angular velocity is measured directly,
\begin{align}
\hat{\omega}_{k+l}^{m}= \tilde\omega_{k+l}.\label{eqn:meaw}
\end{align}

\textit{Uncertainty matrix:} We can represent the actual state at
the $(k+l)$th step using the measured center and perturbations as
follows.
\begin{align}
    C_{k+l}&=\hat{C}_{k+l}^{m} \exp{\Big(S(\zeta_{k+l}^{m})\Big)},\label{eqn:Rkpm}\\
    \omega_{k+l}&=\hat{\omega}_{k+l}^{m}+\delta\omega_{k+l}^{m},\label{eqn:Omegakpm}
\end{align}
for $\zeta_{k+l}^{m},\delta\omega_{k+l}^{m}\in\Re^3$. The
uncertainty matrix is obtained by finding an ellipsoid containing
$\zeta_{k+l}^{m},\delta\omega_{k+l}^{m}$.

We determine the attitude indirectly by comparing the known directions
in the reference frame with measurements in the body frame. So, we need to
transform the uncertainties in the direction measurements into the
uncertainties in the rotation matrix by \refeqn{meaR}. Using the measurement
error model defined in \refeqn{bi}, the actual direction matrix to the
known point ${B}_{k+l}$ is given by
\begin{align}
{B}_{k+l} & = \tilde{B}_{k+l} +
\delta\tilde{B}_{k+l},\label{eqn:Bkl}
\end{align}
where $\delta B_{k+l}=\bracket{S(\nu^1)\tilde b^1,S(\nu^2)\tilde
b^2,\cdots,S(\nu^m)\tilde b^m}\in\Re^{3\times m}$.

The actual matrix giving the known directions ${B}_{k+l}$ and the
actual attitude $C_{k+l}$ at the $(k+l)$th step also satisfy
\refeqn{meaRex};
\begin{gather}
C_{k+l}\Tp L_{k+l}- L_{k+l}\Tp{C}_{k+l}=0,\label{eqn:meaRac}
\end{gather}
where $L_{k+l}=E_{k+l} W_{k+l}{B}_{k+l}\Tp\in\Re^{3\times 3}$.
Substituting \refeqn{Rkpm} and \refeqn{Bkl} into \refeqn{meaRac},
and assuming that the size of the measurement error is sufficiently
small, the above equation can be written as
\begin{align*}
\tilde{L}_{k+l}\Tp \hat{C}_{k+l}^m S(\zeta_{k+l}^m)+
S(\zeta_{k+l}^m)\parenth{\hat{C}_{k+l}^m}\Tp \tilde{L}_{k+l} &
=\parenth{\hat{C}_{k+l}^m}\Tp E_{k+l}W_{k+l}\delta{B}_{k+l}\Tp
-\delta{B}_{k+l}W_{k+l}E_{k+l}\Tp\hat{C}_{k+l}^m.
\end{align*}
Using the identity, $S(x)A+A\Tp S(x)=S(\braces{\tr{A}I_{3\times
3}-A}x)$ for $A\in\Re^{3\times 3},x\in\Re^3$, the above equation can
be written in a vector form.
\begin{align*}
\braces{\tr{\parenth{\hat{C}_{k+l}^m}\Tp
\tilde{L}_{k+l}}-\parenth{\hat{C}_{k+l}^m}\Tp \tilde{L}_{k+l}}
\zeta_{k+l}^m =-\sum_{i=1}^m w_i \braces{\tr{\tilde{b}^i (e^i)\Tp
\hat{C}_{k+l}^m}I_{3\times 3}-\tilde{b}^i (e^i)\Tp
\hat{C}_{k+l}^m}\nu^i.
\end{align*}
Then, we obtain
\begin{align}\label{eqn:zetakpm}
\zeta_{k+l}^m &= \sum_{i=1}^m \mathcal{A}_{k+l}^{m,i}\nu^i,
\end{align}
where
\begin{align}
\mathcal{A}_{k+l}^{m,i} & = -\braces{\tr{\parenth{\hat{C}_{k+l}^m}\Tp
\tilde{L}_{k+l}}-\parenth{\hat{C}_{k+l}^m}\Tp \tilde{L}_{k+l}}^{-1}
 w_i \braces{\tr{\tilde{b}^i (e^i)\Tp
\hat{C}_{k+l}^m}I_{3\times 3}-\tilde{b}^i (e^i)\Tp \hat{C}_{k+l}^m}.
\end{align}
This equation expresses the uncertainty in the measured attitude
as a linear combination of the directional measurement errors.

The perturbation of the angular velocity $\delta\omega_{k+l}^{m}$ is
equal to the angular velocity measurement error $\upsilon_{k+l}$,
since we measure the angular velocity directly. Substituting
\refeqn{Omegakpm} into \refeqn{Omega}, we obtain
\begin{align}\label{eqn:delPikpm}
\delta\omega_{k+l}^m = \upsilon_{k+l}.
\end{align}

Define
$x_{k+l}^m=[(\zeta^m_{k+l})\Tp,\,(\delta\omega^m_{k+l})\Tp]\Tp\in\Re^6$.
Using \refeqn{zetakpm} and \refeqn{delPikpm},
\begin{align*}
x_{k+l}^m & = H_1 \sum_{i=1}^m \mathcal{A}_{k+l}^{m,i} \nu^i_{k+l} +
H_2 \upsilon_{k+l},
\end{align*}
where $H_1=[I_{3\times 3},\, 0_{3\times 3}]\Tp,H_2=[0_{3\times 3},\,
I_{3\times 3}]\Tp\in\Re^{6\times 3}$. Now $x_{k+l}^m$ is expressed as
a linear combination of the measurement errors $\nu^i$ and $\upsilon$.
Using \refeqn{Sk} and \refeqn{Tk}, we can show that each term in the
right hand side of the above equation is in the following
uncertainty ellipsoids.
\begin{align*}
H_1 \mathcal{A}_{k+l}^{m,i} \nu^i_{k+l} & \in
\mathcal{E}_{\Re^6}\parenth{0,H_1 \mathcal{A}_{k+l}^{m,i}S^i_{k+l}\parenth{\mathcal{A}_{k+l}^{m,i}}\Tp H_1\Tp},\\
H_2 \upsilon_{k+l} & \in \mathcal{E}_{\Re^6}\parenth{0,H_2 T_{k+l}
 H_2\Tp}.
\end{align*}
Thus, the uncertainty ellipsoid for $x_{k+l}^{m}$ is obtained as the
vector sum of the above uncertainty ellipsoids. The measurement
update procedure is to find a minimal ellipsoid that contains the
vector sum of those uncertainty ellipsoids. Expressions for a
minimal ellipsoid containing the vector sum of multiple ellipsoids
are presented in~\cite{jo:MaNo1996} and~\cite{jo:DuWaPo2001}.
Using the results, we obtain
\begin{align}
P_{k+l}^m & = \braces{\sum_{i=1}^m \sqrt{\tr{H_1
\mathcal{A}_{k+l}^{m,i}S^i_{k+l}\parenth{\mathcal{A}_{k+l}^{m,i}}\Tp H_1\Tp}}
+\sqrt{\tr{H_2 T_{k+l} H_2\Tp}}}\nonumber\\
& \quad \times \braces{\sum_{i=1}^m \frac{H_1
\mathcal{A}_{k+l}^{m,i}S^i_{k+l}\parenth{\mathcal{A}_{k+l}^{m,i}}\Tp H_1\Tp}{\sqrt{\tr{H_1
\mathcal{A}_{k+l}^{m,i}S^i_{k+l}\parenth{\mathcal{A}_{k+l}^{m,i}}\Tp H_1\Tp}}}
+\frac{H_2 T_{k+l} H_2\Tp}{\sqrt{\tr{H_2 T_{k+l}
H_2\Tp}}}}.\label{eqn:Pmkl}
\end{align}

In summary, the measured uncertainty ellipsoid at the $(k+l)$th step
is defined by \refeqn{meaRex}, \refeqn{meaw}, and \refeqn{Pmkl};
\begin{align}\label{eqn:mea}
    (C_{k+l},\omega_{k+l})\in\mathcal{E}(\hat{C}_{k+l}^{m},\hat{\omega}_{k+l}^{m},P_{k+l}^m).
\end{align}

\vspace*{-2mm}
\subsection{Filtering procedure}
\vspace*{-2mm}

The filtering procedure is to find a new uncertainty ellipsoid
compatible with both the predicted uncertainty ellipsoid and the
measured uncertainty ellipsoid. From \refeqn{flow} and \refeqn{mea},
we know that the state at $(k+l)$th step lies in the intersection of
the two uncertainty ellipsoids:
\begin{align}\label{eqn:filter}
    (C_{k+l},\omega_{k+l})\in\mathcal{E}(\hat{C}_{k+l}^{f},\hat{\omega}_{k+l}^{f},P_{k+l}^f)\bigcap
    \mathcal{E}(\hat{C}_{k+l}^{m},\hat{\omega}_{k+l}^{m},P_{k+l}^m).
\end{align}
However, the intersection of two ellipsoids is not generally an ellipsoid.
We find a minimal uncertainty ellipsoid containing this intersection.
We first obtain equivalent uncertainty ellipsoids in $\Re^6$, and
convert them to uncertainty ellipsoids in $\T\SO$.
We omit the subscript $(k+l)$ in this subsection for convenience.

The uncertainty ellipsoid obtained from the measurements,
$\mathcal{E}(\hat{C}^{m},\hat{\omega}^{m},P^m)$, is identified by
its center $(\hat{C}^{m},\hat{\omega}^{m})$, and the uncertainty
ellipsoid in $\Re^6$:
\begin{align}
(\zeta^m,\delta\omega^m)\in\mathcal{E}_{\Re^6}(0_{6\times
1},P^m),\label{eqn:ellRm}
\end{align}
where $S(\zeta^m)=\mathrm{logm} \parenth{\hat{C}^{m,T} C}\in\so$,
$\delta\omega^m=\omega-\hat{\omega}^m\in\Re^3$. Similarly, the
uncertainty ellipsoid obtained from the flow update,
$\mathcal{E}(\hat{C}^{f},\hat{\omega}^{f},P^f)$, is identified by
its center $(\hat{C}^{f},\hat{\omega}^{f})$, and the uncertainty
ellipsoid in $\Re^6$.
\begin{align}
(\zeta^f,\delta\omega^f)\in\mathcal{E}_{\Re^6}(0_{6\times
1},P^f),\label{eqn:ellRf}
\end{align}
where $S(\zeta^f)=\mathrm{logm} \parenth{\hat{C}^{f,T} C}\in\so$,
$\delta\omega^f=\omega-\hat{\omega}^f\in\Re^3$. Equivalently, an
element
$(C^f,\omega^f)\in\mathcal{E}(\hat{C}^{f},\hat{\omega}^{f},P^f)$, is
given by
\begin{align}
C^f & = \hat{C}^{f} \exp{\Big(S(\zeta^f)\Big)},\label{eqn:Rf}\\
\omega^f & = \hat{\omega}^{f} + \delta{\omega}^f,\label{eqn:Pif}
\end{align}
for some $(\zeta^f,\delta\omega^f)\in\mathcal{E}_{\Re^6}(0_{6\times
1},P^f)$. We find an equivalent expression for \refeqn{ellRf} based
on the center $(\hat{C}^{m},\hat{\omega}^{m})$ obtained from the
measurements.

Define $\hat\zeta^{mf},\delta\hat\omega^{mf}\in\Re^3$ such that
\begin{align}
\hat{C}^f&=\hat{C}^m \exp{\Big(S(\hat\zeta^{mf})\Big)},\label{eqn:errz}\\
\hat{\omega}^f& =\hat{\omega}^m
+\delta\hat\omega^{mf}.\label{eqn:errPi}
\end{align}
Thus, $\hat\zeta^{mf},\delta\hat\omega^{mf}$ represent the
difference between the centers of the two ellipsoids.

Substituting \refeqn{errz}, \refeqn{errPi} into \refeqn{Rf}, \refeqn{Pif}, we obtain
\begin{align}
C^f & = \hat{C}^m \exp{\Big(S(\hat\zeta^{mf})\Big)}
\exp{\Big(S(\zeta^f)\Big)},\nonumber\\
& \simeq \hat{C}^m \exp{\Big(S(\hat\zeta^{mf}+\zeta^f)\Big)},\\
\omega^f & = \hat{\omega}^m + \parenth{\delta\hat\omega^{mf} +
\delta{\omega}^f},
\end{align}
where we assumed that $\hat\zeta^{mf}, \zeta^f$ are sufficiently
small to obtain the second equality. Thus, the uncertainty ellipsoid
obtained by the flow update,
$\mathcal{E}(\hat{C}^{f},\hat{\omega}^{f},P^f)$, is described by the
center $(\hat{C}^m,\hat{\omega}^m)$ obtained from the measurement
and the following uncertainty ellipsoid in $\Re^6$:
\begin{align}
\mathcal{E}_{\Re^6}( \hat{x}^{mf} ,P^f),\label{eqn:ellRmf}
\end{align}
where
$\hat{x}^{mf}=[(\hat\zeta^{mf})\Tp,(\delta\hat\omega^{mf})\Tp]\Tp\in\Re^6$.

We seek a minimal ellipsoid that contains the intersection of two
uncertainty ellipsoids in $\Re^6$.
\begin{align}
\mathcal{E}_{\Re^6}(0_{6\times 1},P^m)\bigcap
\mathcal{E}_{\Re^6}(\hat{x}^{mf} ,P^f)\subset\mathcal{E}_{\Re^6}(\hat{x},P),
\end{align}
where $\hat{x}=[\hat\zeta\Tp,\delta\hat\omega\Tp]\Tp\in\Re^6$.
Using the result presented in~\cite{jo:MaNo1996}, $\hat{x}$ and $P$
can be written as
\begin{align}
\hat{x}=L\hat{x}^{mf},\;\ P=\beta(q) (I-L)P^m, \label{xhatP}
\end{align}
where
\begin{align}
\beta(q) & = 1 + q - (\hat{x}^{mf})\Tp (P^{m})^{-1} L \hat{x}^{mf},
\label{betaq} \\
L & = P^{m} (P^{m} + q^{-1} P^f)^{-1}. \label{LP}
\end{align}
The constant $q$ is chosen such that $\tr{P}$ is minimized. We
convert $\hat{x}$ to points in $\T\SO$ using the common center
$(\hat{C}^{m},\hat{\omega}^{m})$.

In summary, the attitude estimation filter algorithm is given by the
following statement.
\begin{proposition}
The attitude and angular velocity estimates and the new uncertainty
ellipsoid at the $(k+l)$th step are given by
\begin{align}
\hat{C}_{k+l}&=\hat{C}_{k+l}^m \exp{\big(S(\hat\zeta)\big)}, \la{nattes} \\
\hat{\omega}_{k+l}& = \hat{\omega}_{k+l}^m + \delta\hat\omega, \la{nveles} \\
P_{k+l}& = P \la{nunmat},
\end{align}
where $\hat{x}=[\hat\zeta\Tp,\delta\hat\omega\Tp]\Tp\in\Re^6$ and $P\in
\Re^{6\times 6}$ are given by equations (\ref{xhatP})-(\ref{LP}). The
actual state lies in the ellipsoid
\begin{align}
    (C_{k+l},\omega_{k+l})\in\mathcal{E}(\hat{C}_{k+l},\hat{\omega}_{k+l},P_{k+l}),
\end{align}
centered at the estimated attitude and angular velocity states.
\la{filtatt1}
\end{proposition}

The center of the new uncertainty ellipsoid is the estimated state,
considered as point estimates of the attitude and the angular velocity
at the $(k+l)$th step. The uncertainty matrix represents the bounds on
the uncertainty of the estimated state. The size of the uncertainty
matrix characterizes the accuracy of the estimate. If the size of the
uncertainty ellipsoid is small, we conclude that the estimation is
accurate. This estimation is optimal in the sense that the size of the
new uncertainty ellipsoid is minimized. The uncertainty matrix can
also be used to predict the distribution of the uncertainty. The
eigenvector of the uncertainty matrix corresponding to the maximum
eigenvalue shows the direction of the maximum uncertainty.

\vspace*{-2mm}
\subsection{Convergence of Filter}
\vspace*{-2mm}

We now present a sufficient condition under which this estimation
algorithm converges, i.e., the size of the uncertainty matrix
decreases monotonically with measurements. The trace of the
positive definite uncertainty matrix $P$ is the measure of
size used in this analysis.

\begin{theorem}
Let $\Lambda=\mathrm{diag}[\lambda_i]$, $i=1,\ldots,6$, be the
matrix of real positive eigenvalues of $(P^m)^{-1}P^f$ and let
$\lambda_{\rm min}$ be the smallest of the $\lambda_i$. The
estimation algorithm given by Proposition \ref{filtatt1} is
convergent if there exists a constant $c\in(0,1)$ such that the
following inequality holds for every measurement,\be
\|A^f\| < \sqrt{\frac{c(q+\lambda_{\rm min})}{6\chi (P^m)(1+q)}},
\la{sufcond} \ee
where $\|\cdot\|$ denotes the Frobenius norm, $A^f$ is the linear
discrete flow from the previous to the current measurement instant,
and
\[ \chi (P^m)= \sqrt{6+30\kappa (P^m)}, \]
$\kappa (P^m)$ is the condition number of $P^m$.
\label{convthm}
\end{theorem}
\noindent {\em Proof}: For convergence, it is sufficient that the
filtering process is a contraction mapping, which is to say that $\tr{P}<
c\tr{P_0}$ where $c\in (0,1)$ and $P_0$ denotes the uncertainty
matrix of the filtered estimate at the previous measurement instant.
The logic of this proof is as follows.
We first obtain a real-valued function $G(P^m,P^f,q)$, such that
\[ \tr{P} < \|A^f\|^2 G(P^m,P^f,q)\tr{P_0}, \]
where $\|A\|= \tr{\big(A\Tp A\big)}^{\frac12}$ denotes the Frobenius
norm. Our sufficient condition for convergence can then be stated as,
\[\|A^f\| <
\sqrt{\frac{c}{G(P^m,P^f,q)}}. \]
This implies that
\[ \|A^f\|^2 G(P^m,P^f,q)\tr{P_0} < c \tr{P_0}, \]
which is to say that the filter is a contraction mapping.

Using the matrix inversion lemma, we can express the uncertainty matrix
$P$ given by (\ref{xhatP})-(\ref{LP}) as
\[ P= \beta (q)\Big( q(P^f)^{-1}+(P^m)^{-1}\Big)^{-1}. \]
Now we have
\[ \tr {P}= \beta (q)\tr{\Big(q(P^f)^{-1}+(P^m)^{-1}\Big)^{-1}}. \]
>From equations (\ref{betaq})-(\ref{LP}), we have
\[ \beta (q)\le 1+q. \]
>From Lemma \ref{diaglem} in Appendix \ref{diagres}, we know that
$(P^m)^{-1}P^f$ can be diagonalized as
\[ (P^m)^{-1}P^f= U\Lambda U^{-1}, \]
where $\Lambda=\mathrm{diag}[\lambda_i]$, $i=1,\ldots,6$ is a
positive diagonal matrix, i.e., $\lambda_i\in\Re^+$.
Then the uncertainty matrix $P$ is given by
\begin{align}
P &= \beta (q)P^f\Big(qI+(P^m)^{-1}P^f\Big)^{-1} \nn \\
&= \beta (q)P^f U\Big(qI+\Lambda\Big)^{-1}U^{-1}. \la{Pexpd}
\end{align}

We use the above expressions for the uncertainty matrix $P$ and
$\beta(q)$, the Cauchy-Schwartz inequality for the Frobenius norm
of matrices \cite{bo:HJ1991}, and the fact that $P^f$ is a symmetric
positive definite matrix, to obtain:
\begin{align*}
\tr{P} &\le  (1+q)\tr{P^f U(qI+\Lambda)^{-1}U^{-1}} \\
&\le  (1+q)\tr{(P^f)^2}^{\frac12}\| U(qI+\Lambda)^{-1}U^{-1} \| \\
& \le  (1+q)\tr{P^f}\| U\| \|(qI+\Lambda)^{-1}\| \|U^{-1}\| \\
& <  (1+q)\tr{P^f}\| U\| \|U^{-1}\| \tr{(qI+\Lambda)^{-1}}.
\end{align*}
The third step is due to the Frobenius norm being submultiplicative,
as it is equal to the H\"{o}lder 2-norm \cite{bo:Bern}. From the
proof of Lemma \ref{diaglem}, we see that $U=Q\Tp (P^m)^{\frac12}$
where $Q$ is an orthogonal matrix. Therefore, we have
\[
\| U \| = \tr{U\Tp U}^{\frac12}= \tr{P^m}^{\frac12}, \;\
\| U^{-1}\| = \tr{(U\Tp)^{-1}U^{-1}}^{\frac12} = \tr{(P^m)^{-1}}^\frac12.
\]
Let $\sigma_1>\cdots >\sigma_6 >0$ denote the eigenvalues of $P^m$
in descending order. Then
\begin{align*}
\tr{P^m}\tr{(P^m)^{-1}} &= (\sigma_1+\sigma_2+\ldots+\sigma_6)\left(
\frac1{\sigma_1}+\frac1{\sigma_2}+\ldots+\frac1{\sigma_6}\right) \\
&\le  6+6(5)\frac{\sigma_1}{\sigma_6} =6+ 30\kappa (P^m)= \chi^2 (P^m),
\end{align*}
where $\kappa (P^m)=\frac{\sigma_1}{\sigma_6}$ is the condition number of
$P^m$. Therefore, we now have
\begin{align*}
\tr{P}
&\le  (1+q)\tr{(A^f)\Tp A^f P_0}\chi (P^m)\tr{(qI+\Lambda)^{-1}} \\
&\le  (1+q)\tr{\big((A^f)\Tp A^f\big)^2}^{\frac12}\tr{P_0^2}^{\frac12}
\chi (P^m)\tr{(qI+\Lambda)^{-1}} \\
& <  (1+q)\|A^f\|^2\, \chi (P^m)\tr{(qI+\Lambda)^{-1}}\tr{P_0}.\\
& \le  \frac{6 \chi(P^m) (1+q)}{q+\lambda_{\rm min}} \| A^f\|^2 \tr{P_0}
\end{align*}
where, in the last inequality, we used the the fact that $\tr{A^{-1}}\le
\frac{6}{\sigma_{\rm min}(A)}$ when $A\in\Re^{6\times 6}$ is positive
definite and $\sigma_{\rm min}(A)$ is the minimum eigenvalue of $A$.
Using the above inequality, we can let
\[G(P^m,P^f,q)=\frac{6 \chi(P^m) (1+q)}{q+\lambda_{\rm min}}.\]
Our initial discussion then implies that
\[ \|A^f\| < \sqrt{\frac{c}{G(P^m,P^f,q)}} = \sqrt{\frac{c(q+\lambda_{\rm min})}{6\chi (P^m)(1+q)}} \]
is a sufficient condition for convergence.
\qed \\

The rate of convergence is determined by the contraction constant $c\in(0,1)$ for which the inequality (\ref{sufcond}) is satisfied for all measurements. Note that the bound depends on $\lambda_{\rm min}=$minimum eigenvalue of
$(P^m)^{-1}P^f$, the ratio of the relative sizes of the uncertainty
matrices due to the flow and the measurements. Since $P^f=A^f P_0
(A^f)\Tp$, a lower bound for $\lambda_{\rm min}$ can be obtained in terms of $\|A^f\|$ which when combined with the relation (\ref{sufcond}) yields an implicit bound on $A^f$. This implicit bound is a sufficient condition for (\ref{sufcond}) to be satisfied.

The sufficient
condition (\ref{sufcond}) also depends on the condition number of $P^m$,
which suggests that the ease of convergence of this scheme is increased
if $\kappa (P^m)$ is small, i.e., the measurement uncertainty bound is
more spherical and less oblate.
Also note that the size (norm) restriction on $A^f$ imposed by
(\ref{sufcond}) becomes more stringent when $\lambda_{\rm min}$ becomes
smaller as size of the measurement uncertainty (given by $P^m$) becomes
larger. Thus, smaller measurement uncertainty also leads to easier
convergence for the filter, as would be intuitively expected. Other
sufficient conditions for convergence of the filter algorithm can be
obtained using similar analysis.

\section{Conclusion}
\vspace*{-3mm}
A deterministic estimation scheme for the attitude dynamics of a rigid
body in an attitude dependent potential field is presented, with an
assumption of bounded measurement errors. The properties of the proposed
attitude estimation scheme are as follows. This attitude estimator has
no singularities since the attitude is represented by a rotation matrix,
and the structure of the rotation matrix is preserved since it is
updated by group operations in $\SO$ using a Lie group variational
integrator. The proposed attitude estimator is robust to the distribution
of the uncertainty in initial conditions and the measurement noise,
since it is a deterministic scheme based on knowledge of the bounds
in these uncertainties. A sufficient condition for convergence of
this filter has been obtained.

\appendix

\section{Product of two symmetric positive definite matrices}
\la{diagres}
\begin{lemma}
Suppose $A, B\in\Re^{n\times n}$ are symmetric positive definite.
Then, there exists a nonsingular matrix $V\in\Re^{n\times n}$ and a
diagonal matrix $\Lambda=\mathrm{diag}[\lambda_i]$, $i=1,\ldots,n$
such that $\lambda_i\in \Re^+$ and
\be
AB= V^{-1}\Lambda V. \la{diagAB}
\ee
\la{diaglem}
\end{lemma}
\noindent {\em Proof}: Since $A$ is symmetric positive definite, it
can be diagonalized by a real orthogonal matrix $Q_A$ such that
$A=Q_A\Lambda_A Q_A\Tp$ where $\Lambda_A$ is diagonal with positive
real diagonal elements. The symmetric positive definite matrix
square root of $A$ is $A^{\frac12}= Q_A \Lambda_A^{\frac12} Q_A\Tp$.
Since the matrix $A^{\frac12}BA^{\frac12}$ is also symmetric and
positive definite, it is diagonalized by a real orthogonal matrix
$Q$ such that \be A^{\frac12}BA^{\frac12}= Q\Lambda Q\Tp,
\la{ABAdiag} \ee where $\Lambda$ is diagonal with positive real
diagonal elements. Define a non-singular matrix $V= Q\Tp
A^{-\frac12}\in\Re^{n\times n}$. Now consider
\begin{align*}
V A B V^{-1}&= Q\Tp A^{1/2} B A^{1/2} Q= Q\Tp Q \Lambda Q\Tp Q= \Lambda,
\end{align*}
which yields (\ref{diagAB}). \qed \\
An alternate proof of the above statement is given in
\cite{bo:Bern}, while a more general result applicable to
a product of two Hermitian matrices is given in \cite{bo:HJ1991}.

\bibliographystyle{elsart-num}
\bibliography{glates}

\end{document}